\documentstyle{amsppt}
\nologo
%
%
\font\b=cmr10 scaled \magstep4
\def\bigzerou{\smash{\kern-25pt\lower1.7ex\hbox{\b 0}}}
\hsize=360pt
\vsize=500pt
\hbadness=5000
\tolerance=1000
\NoRunningHeads

\def\Q{\Bbb Q}

\def\F5{\Bbb F_5}

\def\no{\noindent}

\def\Ga{\Gamma}

\def\part{\partial}

\magnification=\magstep1
\topmatter
\title A rational map between two threefolds \endtitle
\author
Kenichiro Kimura
\endauthor

\endtopmatter
\vskip -3ex
\noindent\hskip 2em

\vskip 3ex
\CenteredTagsOnSplits
%
\document

Let $V_{33}\subset \Bbb P^5$ be the threefold defined by the following 
2 cubic equations:

$$ \align X_0^3+X_1^3+&X_2^3+X_3^3=0 \\
          &X_2^3+X_3^3+X_4^3+X_5^3=0 .\endalign $$
Let $\widetilde{V}_{33}$ be the blow up of $V_{33}$ along its singular
locus. In \cite{12} it is shown that the third betti number of  
$\widetilde{V}_{33}$ is 2 and the $L$-series $L(H^3( \widetilde{V}_{33}),s)=
L(f,s)$ with $f$ the newform of weight 4 on $\Gamma_0(9)$. This is also
the $L$-series of a Hecke character of the field $\Q(\sqrt{-3})$. 

\no Let $E\subset \Bbb P^2$ be the curve defined by

$$ X_0^3+X_1^3+X_2^3=0.$$ In Remark 4.5 
of \cite{4} it is mentioned that 
there is a piece of the cohomology $H^3(E^3)$ which 
has the same $L$-series as that of $H^3( \widetilde{V}_{33})$. 
According to Tate's
conjecture there should be a correspondence between $V_{33}$ and $E^3$.
We give such a correspondence explicitly. Other relatives of these
varieties are the fibered squares of the universal curves for 
$\Ga(3)$ and $\Ga_0(9)$. In \cite{8} and in the section 
13 of \cite{9} Schoen constructs correspondences 
between these varieties and $E^3$. 

\no We refer the reader to Remark 4.5 
of \cite{4} for a list of threefolds which have a two dimensional Galois
representation in $H^3$ whose $L$-series is associated to an elliptic
modular form of weight 4, and known correspondences among those varieties.
We include the relevant works in the reference.

\proclaim{Theorem} There is a dominant rational map from $E^3$ to $V_{33}$
of degree 3. \endproclaim
\demo{Proof} Consider the affine piece $\{X_3\neq 0\}$ of $V_{33}$
which is equal to
$$ \text{Spec}\Q[X_0,X_1,X_2,X_4,X_5]/
( X_0^3+X_1^3+X_2^3+1,\,X_2^3+1+X_4^3+X_5^3). $$
On the other hand, the affine open $(E-\{X_2=0\})^3$ of 
$E^3$ is given by
$$\text{Spec}\Q[x_1,y_1,x_2,y_2,x_3,y_3]/(x_1^3+y_1^3+1,x_2^3+y_2^3+1,
x_3^3+y_3^3+1).$$  
Then there is a map from $(E-\{X_2=0\})^3$
to $V_{33}$ induced by the ring homomorphism
$$ X_0\mapsto -x_1y_3,\,X_1\mapsto -y_1y_3,\,X_2\mapsto x_3,\,
X_4\mapsto -x_2y_3,\,X_5\mapsto -y_2y_3.$$ We denote by $\Q(E^3)$
and $\Q(V_{33})$ the function fields of $E^3$ and of $V_{33}$
respectively. Then one can see that the field $\Q(E^3)$ is generated
over $\Q(V_{33})$ by $y_3$ with the equation $y_3^3+X_2^3+1=0$.
So this map is dominant of degree 3. \enddemo

\subsubhead{acknowledgment}\endsubsubhead
 The author is grateful to the referee for
comments on earlier version of this paper which improved the exposition.

$\quad$

Institute of Mathematics

University of Tsukuba

Tsukuba

305-8571

Japan 

kimurak\@math.tsukuba.ac.jp  


\Refs
\widestnumber\key{10}














\ref
\key 1
\paper An exceptional isomorphism between modular varieties
\by Ekedahl,T.; van Geemen, B. 
\jour Progr. Math. 
\yr 1991
\pages 51-74
\vol 89
\endref


\ref
\key 2
\paper On the Geometry and Arithmetic of Some Siegel Modular Threefolds
\by  van Geemen, B.; Nygaard, N.
\jour J. of Number Theory
\yr 1995
\pages  45-87
\vol 53
\endref


\ref
\key 3
\paper An Isogeny of $K3$ surfaces
\jour mathAG/0309272
\yr 2003
\by van Geemen, V.; Top, J.
\endref

\ref
\key 4
\paper The modularity of the Barth-Nieto quintic and its relatives
\by   Hulek, K.; Spandaw, J.; van Geemen, B.; van Straten, D.
\jour Adv. Geom.
\yr 2001
\pages 263-289 also available as AG/0010049
\vol 1
\endref

\ref
\key 5
\paper Cubic exponential sums and Galois representations
\by  Livne, R.
\jour Contemp. Math.
\yr 1987
\pages  247-261
\vol 67
\endref

\ref
\key 6
\paper The modularity conjecture for rigid Calabi-Yau threefolds over $\Q$
\jour J. Math. Kyoto Univ.
\vol 41
\yr 2001
\pages 403-419
\by Saito, M-H.; Yui, N.
\endref


\ref
\key 7
\paper On the geometry of a special determinantial hypersurface associated
to the Mumford-Horrocks vector bundle
\by  Schoen, C.
\jour J. Reine Angew. Math.
\yr 1986
\pages 85-111
\vol 364
\endref


\ref
\key 8
\paper Zero cycles modulo rational equivalence for some varieties over fields of transcendence degree one
\by  Schoen, C.
\jour Proc. Sympos. Pure Math.
\yr 1987
\pages 463-473
\vol 46, Part 2
\endref


\ref
\key 9
\paper On the computation of the cycle class map for nullhomologous cycles over the algebraic closure of a finite field
\by  Schoen, C.
\jour Ann. Sci. Ecole Norm. Sup. (4)
\yr 1995
\pages 1-50
\vol 28
\endref

\ref
\key 10
\paper A quintic hypersurface in $\Bbb P^4$ with 130 nodes
\jour Topology
\yr 1993
\vol 32
\pages 857-864
\by van Straten, D.
\endref

\ref
\key 11
\paper The $L$-series of certain rigid Calabi-Yau threefolds
\jour J. Number Theory
\vol 81
\yr 2000
\pages 509-542
\by Verrill, H.
\endref


\ref
\key 12
\paper New examples of threefolds with $c_1=0$
\jour Math. Z.
\yr 1990
\vol 203
\pages 211-225
\by Werner, J.; van Geemen, B.
\endref

\ref
\key 13
\paper The arithmetic of certain Calabi-Yau varieties over number fields
\jour NATO Sci. Ser. C Math. Phys. Sci.
\vol 548
\yr 2000
\pages 515-560
\by Yui, N.
\endref




\endRefs
\enddocument